\documentclass[letterpaper,draftcls,onecolumn,twoside,11pt]{IEEEtran}
\usepackage{amsxtra}
\usepackage{amssymb}
\usepackage{amsmath}
\usepackage[dvips]{graphicx} 
\usepackage[bf]{subfigure}
\usepackage{cite}

\def\tr{\mathop{\mathrm{tr}}}
\def\b{\boldsymbol{b}}
\def\e{\boldsymbol{e}}
\def\w{\boldsymbol{w}}
\def\x{\boldsymbol{x}}
\def\y{\boldsymbol{y}}
\def\xhat{\boldsymbol{\widehat{{x}}}}

\title{Beyond Thresholding:
       Analysis and Improvements for Deterministic Parameter Estimation}

\author{Baris I. Erkmen, \IEEEmembership{Student Member, IEEE} and
        Vivek K. Goyal, \IEEEmembership{Senior Member, IEEE}
\thanks{This work was supported in part by the
        Texas Instruments Leadership University Consortium Program.}
\thanks{B. I. Erkmen (email: erkmen@mit.edu) and V. K. Goyal (email: vgoyal@mit.edu) are with the Department of Electrical Engineering and Computer Science and
        the Research Laboratory of Electronics, Massachusetts Institute of Technology, 77 Massachusetts Avenue, Cambridge, MA 02139, USA.}
}

\begin{document}

\maketitle

\begin{abstract}
Hard-threshold estimators are popular in
signal processing applications. We provide a detailed study of using
hard-threshold estimators for estimating an unknown deterministic
signal when additive white Gaussian noise corrupts observations.
The analysis, depending heavily on Cram{\'e}r-Rao bounds, motivates
piecewise-linear estimation as a simple improvement to hard thresholding.
We compare the performance of two piecewise-linear estimators to a
hard-threshold estimator.  When either piecewise-linear estimator is
optimized for the decay rate of the basis coefficients, its performance
is better than the best possible with hard thresholding.
\end{abstract}


\begin{center}
EDICS Categories:  SSP-PARE, SSP-REST, SSP-PERF
\end{center}

\newpage

\section{Introduction} 
\IEEEPARstart{R}{emoving} noise from signals (``denoising'')
is a problem central to many engineering disciplines.
As summarized nicely by Moulin and Liu~\cite{MoulinL:99},
most methods fall into at least one of three categories:
Bayesian techniques, which assume a probabilistic prior for the
unknown signal and minimize an error measure given the observations;
minimax techniques, which are designed for good worst-case performance over
some broad class of signals; and
techniques based on the minimum description length (MDL) principle.
In the electrical engineering literature,
the Bayesian and MDL approaches are more common
than the minimax approach familiar to statisticians.

In many fields---especially image processing and geophysics---the
use of the wavelet domain is prominent,
regardless of which of the three approaches is undertaken.
For Bayesian techniques, the wavelet domain is convenient because it allows
low-complexity diagonal~\cite{SimoncelliA:96,ChipmanKM:97,PortillaSWS:03}
or nearly-diagonal~\cite{SendurS:02_SP} estimators to be used with
little loss in performance.
Minimax estimation performance is intimately connected with
nonlinear approximation~\cite{Donoho:93},
so the approximation power of wavelet bases for many classes of signals
make them appropriate~\cite{AntoniniBMD:92,DeVoreJL:92,DeVore:98}.
Finally, the success of wavelet-based
compression makes wavelet representations suitable for generation of
regularization terms in MDL~\cite{ChangYV:00b,HansenY:00}.
All three approaches, under appropriate conditions,
justify simple hard threshold or soft threshold (shrinkage) estimators.

In this paper we study the classical problem of
estimating a signal in the presence of additive white Gaussian noise (AWGN)
with the goal of minimizing mean-squared error (MSE)\@.
Rather than applying the Bayesian formulation,
we cast this as the estimation of a non-random parameter vector.
This allows us to explain the performance of hard threshold
estimators through bias-variance trade-off and the Cram{\'e}r-Rao bound (CRB).
The (biased) CRB provides more insight than the standard ``oracle'' bound
of~\cite[Ch.~10]{Mallat:99}.  Shaping the bias and resulting MSE inspires
the analysis and optimization of two alternatives to hard thresholding.
Both are piecewise-linear functions, and one has been proposed
previously as the ``semisoft shrinkage'' estimator~\cite{GaoB:95}.
Our focus is not on the invention of simple estimators, but rather on the fact that the performance of such estimators can be better understood with the analysis presented herein. Furthermore, we show that the degrees of freedom in these estimators can be optimized, given the decay rate of coefficients, to achieve lower average estimation error than that incurred by hard thresholding.

The paper is organized as follows.
In Section~\ref{sec:background} we review the definition of bias
and estimation error bounds for both unbiased and biased estimators.
We then analyze hard-threshold estimators---notably explaining their
performance using Cram{\'e}r-Rao Bounds (CRBs)---in Section~\ref{sec:ht}.
Inspired by this way of understanding hard-threshold estimators,
we analyze two alternative estimators in Section~\ref{sec:alternatives}.
In particular, we show that these estimators can be optimized for
the decay rate of the unknown deterministic parameter vector,
resulting in uniform improvement over hard thresholding.
We also discuss the limiting cases that relate the alternative
estimators to hard-threshold estimators.
Finally, Section~\ref{sec:discussion}
provides a discussion of the key results that emerge from our analysis.

\section{Background on Estimation Error Bounds}
\label{sec:background}
\subsection{Estimation in White Gaussian Noise}
\label{sec:AWGNestimation}

In this paper we consider non-random signal estimation when the observation
is the signal plus white Gaussian noise.
In particular, assume that the observed signal is expanded in some basis
of our choice, and the $N \in \mathbb{Z}^{+}$ basis coefficients
of interest are stacked in a $N \times 1$ vector $\y \in \mathbb{R}^{N}$, such that
\begin{equation}
\y = \x  + \w ,
\end{equation}
where $\x\in \mathbb{R}^{N}$ and $\w \in \mathbb{R}^{N}$ are $N \times 1$ column vectors
representing the corresponding signal and noise basis coefficients respectively,
when expanded in the same basis.\footnote{Throughout this paper we are going
to assume $N$ is finite, although it can be arbitrarily large.}
Here, $\x$ is deterministic yet unknown, and $\w$
is a random, zero-mean ($E[\w] = 0 $)
white Gaussian noise vector with correlation matrix
$E[\w \w^{T}] = \sigma_{w}^{2} I_N$,
where $I_N$ is the $N \times N$ identity matrix.
Thus $\y$ is a Gaussian random vector with mean $\x$
and covariance matrix $\sigma_{w}^2 I_N$;
i.e. the probability density function of $\y$ is
\begin{equation}
p(\y ; \x) \sim \mathcal{N}(\x,\sigma_{w}^{2} {I_{N}})
 = \frac{1}{(2 \pi \sigma_{w}^{2})^{N/2}} e^{- \lVert \y-\x\rVert^{2}/(2\sigma_{w}^2)} .
\end{equation}

Before embarking on our analysis let us review and establish some notation.
An estimator for $\x$, denoted by
$\xhat(\y)$,
is a deterministic function of $\y$
that maps an observation vector in $\mathbb{R}^{N}$
into the parameter space $\mathbb{R}^{N}$.\footnote{The
observation space and parameter space need not be identical in general.}
Given an estimator, we define the error as
\begin{equation}
\e(\y) \equiv \xhat(\y) - \x ,
\end{equation}
which is a $N \times 1$ random vector.
Its mean value is termed the bias of the estimator and is denoted with
\begin{equation}
\b(\x) \equiv E[\xhat(\y)] - \x .
\end{equation}
Note that the bias of an estimator is in general a function of the parameters
that are being estimated (i.e.\ $\x$),
so in general it is not trivial to arbitrarily modify or eliminate
the bias of an estimator.

The performance of an estimator is often assessed by its $N \times N$, 
positive-definite error correlation matrix
\begin{equation}\
\Lambda_{e}(\x) \equiv E[\e(\y) \e^{T}(\y)].
\label{errorcor}
\end{equation}
Note that the $\ell_{2}$-norm of the error vector can be obtained
from the error correlation matrix by taking its trace; i.e.
\begin{equation}
\mbox{mse}(\x) \equiv E \left [\lVert \xhat(\y) - \x \rVert^{2} \right ] = \tr \left  (\Lambda_{e} \right ) .
\end{equation}

\subsection{Unbiased Estimators}
\label{sec:unbiased}

An estimator is {\em unbiased} if it satisfies
$\boldsymbol{b}(\boldsymbol{x}) = 0$
for all $\boldsymbol{x}$.
It is well known in estimation theory \cite[Ch.~2]{VanTrees:01a}
that the Cram\'{e}r-Rao bound yields a global performance bound
for all unbiased estimators as
\begin{equation}
\Lambda_{e} - {I}_{\boldsymbol{y}}^{-1}(\boldsymbol{x})\geq 0 , 
\label{CRB:UB}
\end{equation}
where `$\geq$' indicates that the matrix on the left hand side is
positive semi-definite. Here ${I}_{\boldsymbol{y}}(\boldsymbol{x})$
is the $N \times N$ Fisher Information matrix with elements
\begin{equation}
[{I}_{\boldsymbol{y}}(\boldsymbol{x})]_{n,m}
 = -E\left[\tfrac{\partial^{2}}{\partial x_{n} \partial x_{m}}
         \text{ln} \, p(\boldsymbol{y};\boldsymbol{x}) \right ] ,
\label{FI}
\end{equation}
for $n,m = 1,\dots,N$.
The Cram\'{e}r-Rao bound is always a lower bound
(in the positive semi-definite sense),
however it may not be possible to satisfy it with equality.
In particular, the left hand side of \eqref{CRB:UB}
is equal to $0$ if and only if the efficient estimator
\begin{equation}
\xhat(\boldsymbol{y})
 = \boldsymbol{x} + {I}_{\boldsymbol{y}}^{-1}(\boldsymbol{x})
\Bigl (\nabla_{\boldsymbol{x}} \ln\, p(\boldsymbol{y};\boldsymbol{x})\Bigr)^{T} ,
\label{effestimator}
\end{equation}
where $\nabla_{\boldsymbol{x}} \equiv [\partial/\partial x_{1},\, \ldots,\, \partial/\partial x_{N} ]$,
exists; i.e. the right hand-side of \eqref{effestimator}
must be independent of $\boldsymbol{x}$.
This is indeed the case for the AWGN problem studied in this paper.
In particular, we find that
\begin{equation}
{I}_{\boldsymbol{y}}^{-1}(\boldsymbol{x}) = \sigma_{w}^{2} {I}_{N} ,
\label{CRBUB:AWGN}
\end{equation}
and $\Lambda_{e} = {I}_{\boldsymbol{y}}^{-1}(\boldsymbol{x})$ is satisfied when the maximum-likelihood estimator is used, i.e. when
\begin{equation}
\xhat(\boldsymbol{y}) = \boldsymbol{y} .
\label{effestimator:AWGN}
\end{equation}
Equations \eqref{CRBUB:AWGN} and \eqref{effestimator:AWGN} imply that the $\ell_{2}$-norm of the error for any unbiased estimator is no less than $N \sigma_{w}^{2}$ and this minimum is achieved by the maximum-likelihood estimator, which is not only a linear estimator, but also the trivial identity function. Furthermore, when the maximum-likelihood estimator is employed, the estimation errors for each element in $\boldsymbol{x}$ are statistically independent.

For unbiased estimators, the CRB on the error covariance matrix is
\emph{independent} of the basis representation, so in particular
a decomposition into a wavelet basis has no advantage over any other basis~\cite[Ch.~2]{VanTrees:01a}. This picture changes significantly however, when we turn our attention to biased estimators.

\subsection{Biased Estimators}
\label{sec:biased}

The set of all estimators that satisfy $\b(\x) \neq 0$ for some $\x$ constitute the class of biased estimators. The Cram\'{e}r-Rao bound for biased estimators is given by \cite[Ch.~2]{VanTrees:01a}
\begin{align}
\lefteqn{\Lambda_{e} - \boldsymbol{b}(\boldsymbol{x}) \boldsymbol{b}^{T}(\boldsymbol{x})} \nonumber \\ && - \Bigl ({I_{N}} + \nabla_{\! \boldsymbol{x}} \boldsymbol{b}(\boldsymbol{x}) \Bigr) {I}_{\boldsymbol{y}}^{-1} (\boldsymbol{x}) \Bigl ({I_{N}} + \nabla_{\! \boldsymbol{x}} \boldsymbol{b}(\boldsymbol{x}) \Bigr)^{\!T} \geq 0 , \label{CRB:B}
\end{align}
where `$\geq$', once again, denotes the positive semi-definiteness of the matrix on the left hand side. Unfortunately, the Cram\'{e}r-Rao bound in the biased case is less useful than it is in the unbiased case. In particular, \eqref{CRB:B} is interpreted as a lower bound for all estimators with bias $\boldsymbol{b}(\boldsymbol{x})$, but it is not guaranteed that multiple estimators can satisfy a given bias function $\boldsymbol{b}(\boldsymbol{x})$ (unless $\boldsymbol{b}(\boldsymbol{x})$ is trivially a constant). Furthermore, in general the bound cannot be satisfied with equality. However, as we shall see shortly, the Cram\'{e}r-Rao bound can provide valuable insight into the performance of particular estimators. 

It is worth emphasizing that the Cram\'{e}r-Rao bound in \eqref{CRB:B} depends on only three quantities: The Fisher Information matrix, the bias, and the gradient of the bias in the parameter space. The Fisher Information matrix is a function of the probability density function of the observed data and therefore is fixed for a given problem setup; for example the Fisher Information matrix for the AWGN problem analyzed in this paper is $\sigma_{w}^{-2} {I_{N}}$. Thus the Cram\'{e}r-Rao bound can be thought as a function of only the latter two quantities.

Because there is no global lower bound on the performance of biased estimators, a general treatment is not possible. Hence, the utility of the biased Cram\'{e}r-Rao bound is best demonstrated via an example. Due to its popularity in wavelet-based estimation techniques, we shall first focus our attention on hard thresholding estimators.

\section{Analysis of Hard Thresholding}
\label{sec:ht}
\begin{figure}
\centering
\includegraphics[width=2in]{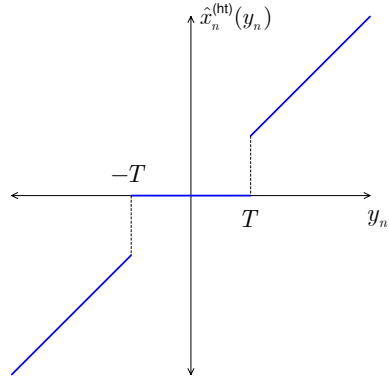}
\caption{The input/output relation of a hard thresholding estimator.}
\label{HT:ESTFIG}
\end{figure}

A hard thresholding estimator acts on each basis coefficient observation $\{ y_{n} \}_{n=1}^N$ independently and estimates the true value of each signal basis coefficient according to
\begin{equation}
\hat{x}_{n}^{\text{(ht)}} (y_{n})
= \begin{cases} 0 , & \text{if } |y_{n}|< T ; \\
y_{n} , & \text{if } |y_{n}|\geq T , \end{cases}
\label{HT:EST}
\end{equation}
where $T > 0$ is called the \emph{threshold}, as shown in Figure~\ref{HT:ESTFIG}. Note that, because the estimator acts on each coefficient separately and because the noise on each coefficient is independent in the AWGN estimation problem, we can restrict our analysis to the scalar estimation case with no loss of generality.

\begin{figure}
\centering
\subfigure[Bias]{
		\includegraphics[width=2.7in]{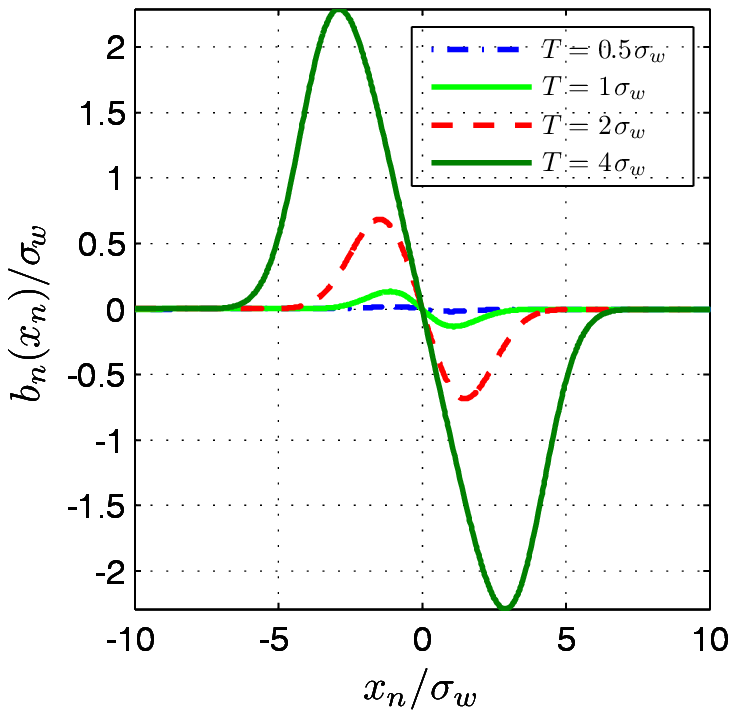}
		\label{HT:BIASFIG}
	}
\subfigure[Mean-square error]{
		\includegraphics[width=2.7in]{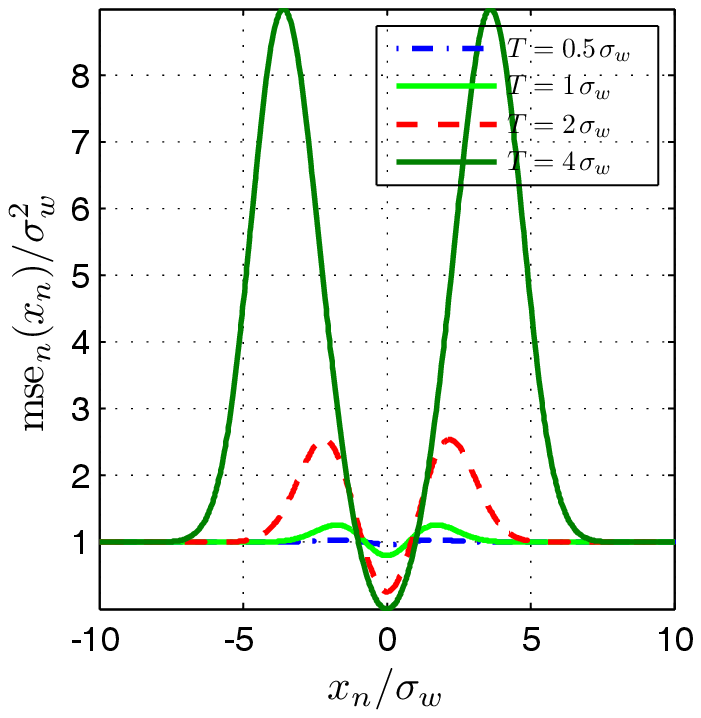}
		\label{HT:MSEFIG}
	}
	\caption{The bias and mean-square error of a hard thresholding estimator as a function of the true value of the parameter. The mean-square error is normalized to the variance of the noise, whereas the bias and the parameter values are normalized to the standard deviation of the noise, such that all axis variables are dimensionless.}
\end{figure}

The hard thresholding estimator differs from the maximum-likelihood estimator in \eqref{effestimator:AWGN} only when the observation $y_{n}$ has absolute value less than the threshold $T$. In this case the hard thresholding estimator estimates the true value of the underlying signal as $0$, and this squelching action introduces a bias given by\footnote{Because the explicit expressions of the bias and mean-square error are cumbersome, we shall defer them to Appendix~\ref{app2} and state here instead concise integral expressions for these quantities.}
\begin{equation}
b_{n}^{\text{(ht)}}(x_{n}) = - \int_{-T}^{T} {y \, p_{w}(y-x_{n}) \, {\rm d}y } ,
\label{HT:BIAS}
\end{equation}
where $p_{w}(y) = \exp(-y^{2}/2 \sigma_{w}^2)/\sqrt{2 \pi \sigma_{w}^2}$. This bias is plotted in Figure~\ref{HT:BIASFIG} for various threshold values. Note that the bias is anti-symmetric; i.e. $b_{n}^{\text{(ht)}}(x_{n}) = -b_{n}^{\text{(ht)}}(-x_{n})$. 

The mean-square error is found to be
\begin{align} 
\text{mse}_{n}^{\text{(ht)}}(x_{n}) = \sigma_{w}^{2}- 2 x_{n} b_{n}^{\text{(ht)}}&(x_{n}) \nonumber \\ -&\int_{-T}^{T} {y^{2} \, p_{w}(y-x_{n}) \, {\rm d}y} .
\label{HT:MSE}
\end{align}
It is worth identifying the terms contributing to this expression. The first term on the right hand side of \eqref{HT:MSE} is the mean-square error obtained when a maximum-likelihood estimator is used. The middle term is always positive and therefore always increases the mean-square error above the minimum mean-square error achievable with an unbiased estimator. On the other hand, the rightmost term is negative, hence it reduces the mean-square error. Therefore, the overall performance of a hard thresholding estimator, relative to the maximum-likelihood estimator, is determined by which of the latter two terms has greater magnitude. Figure~\ref{HT:MSEFIG} plots the mean-square error normalized to the variance of the noise, against the true value of the parameter $x_{n}$ normalized to the standard-deviation of the noise. The figure shows that the performance of a hard thresholding estimator depends on the value of $x_{n}$. Recalling that the maximum-likelihood estimator has mean-square error equal to $\sigma^{2}_{w}$ independent of $x_{n}$, we observe from the plots that when $x_{n}$ is within approximately one standard-deviation of $0$, the mean-square error of the hard thresholding estimator is less than that of the maximum-likelihood estimator. On the other hand, if $x_{n} \gg T$, then the mean-square error approaches that of the maximum-likelihood estimator, because the probability that the noise will push the observation into the regime of thresholding becomes very small. However, in the intermediate regime of $x_{n}$, when the true value of the signal is on the same order as the threshold, the mean-square error of the thresholding estimator is worse than the maximum-likelihood estimator, because the noise can push the observation to either side of the threshold, leading to significant errors in the estimate.

\begin{figure}
\centering
\includegraphics[width=2.7in]{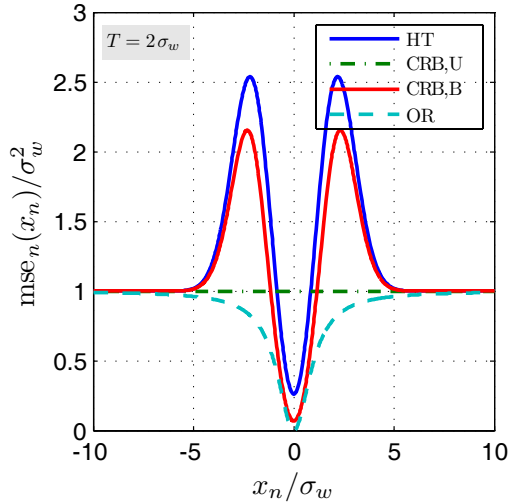}
\caption{The mean-square error of a hard thresholding estimator with $T=2 \sigma_{w}$, together with the unbiased Cram\'{e}r-Rao bound (CRB,U), the biased Cram\'{e}r-Rao bound (CRB,B) for all estimators with bias given by $b_{n}^{\text{(ht)}}(x_{n})$, and the ``optimal oracle" bound (OR) derived in \cite[Ch.~10]{Mallat:99} (see Appendix~\ref{app3}). Because the hard thresholding estimator is biased, its mean-square error can go lower than the unbiased Cram\'{e}r-Rao bound.}
\label{HT:LBFIG}
\end{figure}

We can utilize the Cram\'{e}r-Rao bound for biased estimators to gain further insight into hard thresholding. Figure~\ref{HT:LBFIG} compares the normalized mean-square error of a hard thresholding estimator to the unbiased and biased Cram\'{e}r-Rao bounds, as well as the ``optimal oracle" lower bound obtained in \cite[Ch.~10]{Mallat:99} and reproduced in Appendix~\ref{app3} for convenience.
Notice from the figure that the oracle bound is a very weak lower bound which does not capture the oscillatory behavior of the hard thresholding mean-square error. On the other hand, the Cram\'{e}r-Rao bound---a lower bound for all estimators with bias given by \eqref{HT:BIAS}---follows the same oscillatory trend as the hard thresholding mean-square error. In the scalar AWGN case, the biased Cram\'{e}r-Rao bound simplifies to
\begin{equation} 
\text{mse}_{n}(x_{n}) \geq b_{n}^{2}(x_{n}) + \sigma_{w}^{2} \Bigl (1 + \tfrac{\partial}{\partial x_{n}} b_{n}(x_{n}) \Bigr)^{2} ,
\label{CRB:BSC}
\end{equation}
which depends {\em only} on the bias and its derivative with respect to $x_{n}$. Hence, the oscillatory behavior of the hard thresholding mean-square error is primarily a consequence of its {\em bias}. From Figure~\ref{HT:BIASFIG} we can verify that the improvement in mean-square error for $x_{n} \approx 0$, is due to $\partial b_{n}^{\text{(ht)}}(x_{n}) /\partial x_{n} < 0$, whereas both $|b_{n}(x_{n})|\gg 0$ and $\partial b_{n}^{\text{(ht)}}(x_{n}) /\partial x_{n}  > 0$ contribute to the peak in the mean-square error.

We can further exploit the bias dependence of the mean-square error to improve the performance of the estimator on a sequence of coefficients with a given decay rate. We develop this in the next section.

\section{Alternatives to Hard Thresholding}
\label{sec:alternatives}
As evidenced in hard thresholding estimators, the bias plays a significant role in determining the mean-square error behavior of an estimator. This dependence has been exploited in previous work to improve the mean-square error performance---with respect to the unbiased Cram\'{e}r-Rao bound---over the entire parameter space \cite{Eldar:06}. Here, we shall instead aim at improving the average error performance of the hard thresholding estimator when applied to multiple basis coefficients with known decay rate. 

\begin{figure}
\centering
\subfigure[Piecewise-linear]{
		\includegraphics[width=2in]{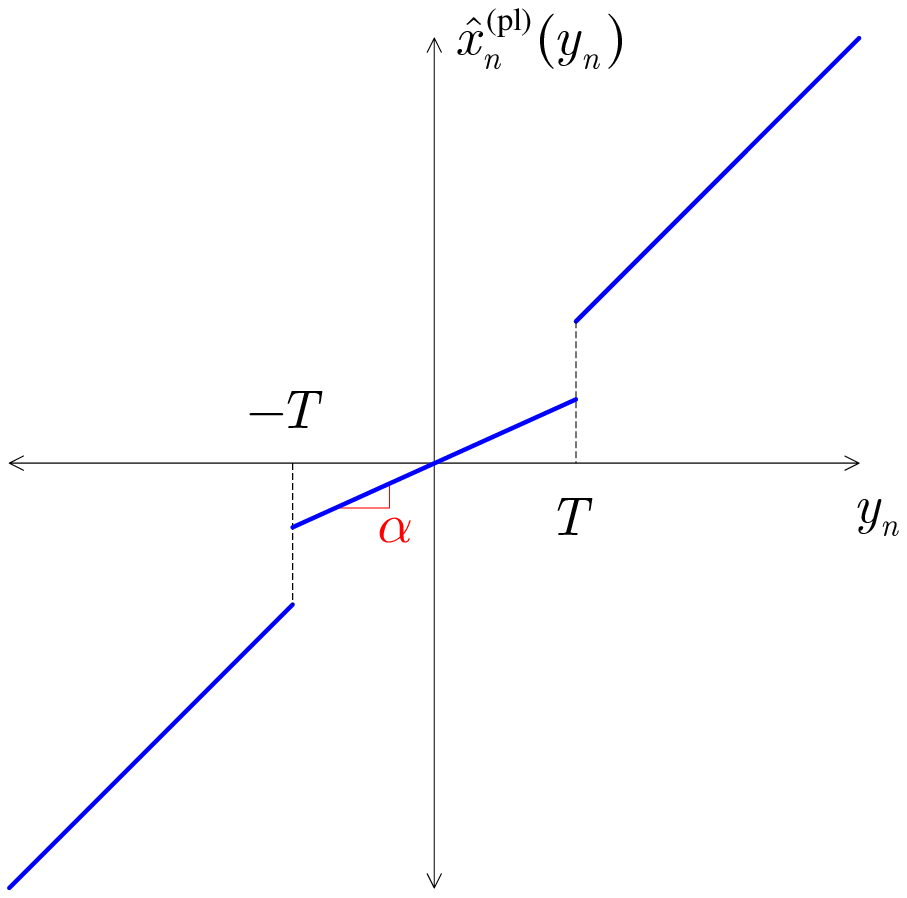}
		\label{PL:ESTFIG}
	}
\subfigure[Semisoft shrinkage]{
		\includegraphics[width=2in]{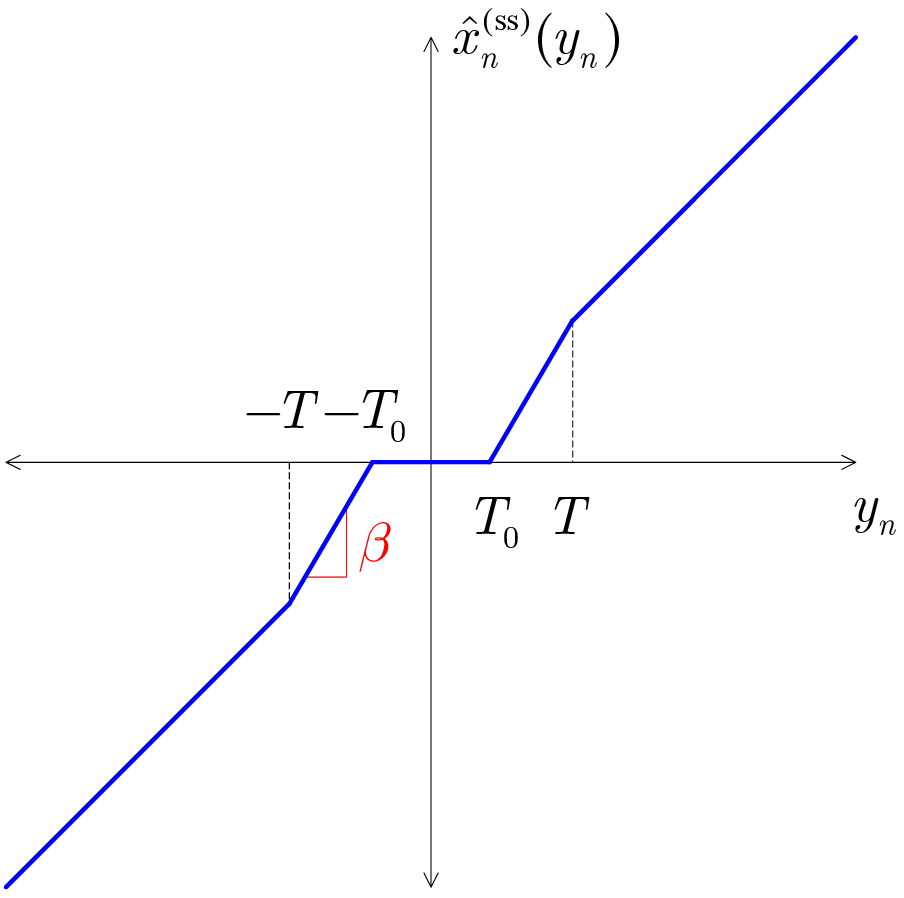}
		\label{SS:ESTFIG}
	}
\caption{The input/output relations for two alternative piecewise-linear estimators.}
\end{figure}

In this section we consider two piecewise-linear estimators that generalize the hard thresholding estimator. The first is given by
\begin{equation}
\hat{x}_{n}^{\text{(pl)}} (y_{n}) = \begin{cases} \alpha y_{n}  , & \text{if } |y_{n}|< T \\ y_{n} , & \text{if } |y_{n}|\geq T \end{cases} , \label{PL:EST}
\end{equation}
where $\alpha \in [0,1]$ and $n = 1,\dots,N$. From the input/output relation of the piecewise-linear estimator shown in Figure~\ref{PL:ESTFIG}, it is clear that $\alpha=0$ corresponds to hard thresholding, whereas $\alpha=1$ yields the maximum-likelihood estimator. Thus, the slope of the line segment over $y_{n} \in [-T, T]$ is a degree of freedom in the piecewise-linear estimator that encompasses both the maximum-likelihood estimator and the hard thresholding estimator as special instances.

\begin{figure}
\centering
\subfigure[Bias]{
		\includegraphics[width=2.7in]{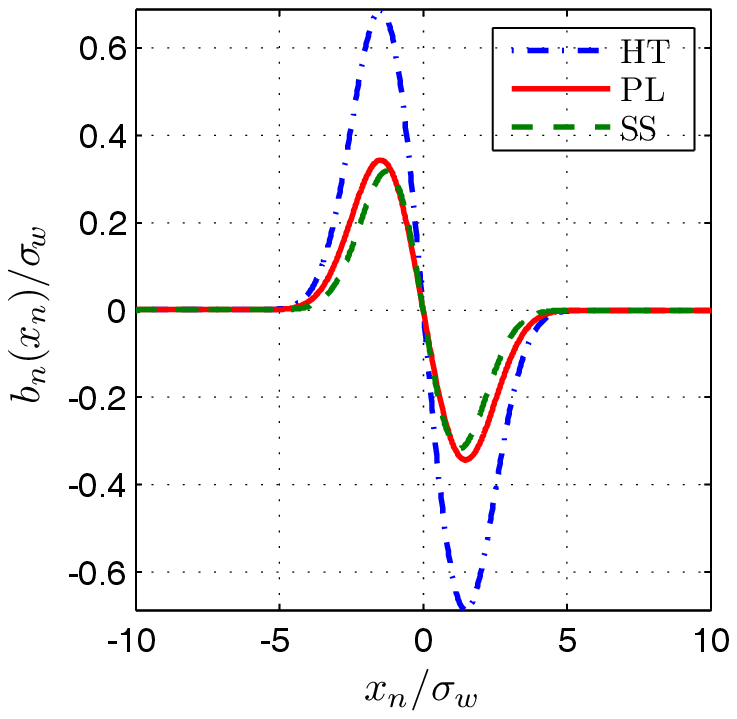}
		\label{PL:BIASFIG}
	}
\subfigure[Mean-square error]{
		\includegraphics[width=2.7in]{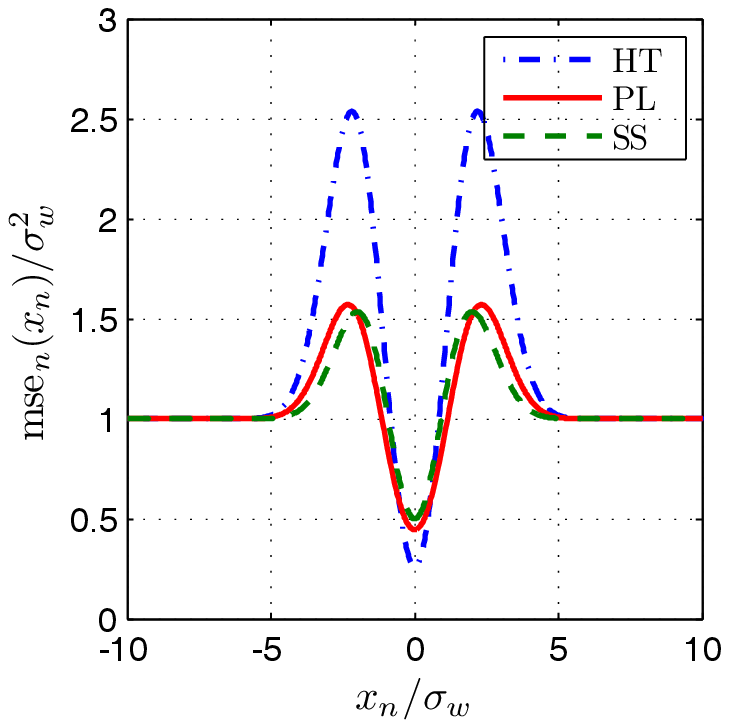}
		\label{PL:MSEFIG}
	}
	\caption{The bias and mean-square error of the two piecewise-linear estimators (PL, SS) and the hard thresholding (HT) estimator are plotted as a function of the true value of the signal coefficient. All axis variables are normalized to be dimensionless. $T=2 \sigma_{w}$ for all of the plots, $\alpha=0.5$ is the slope of the piecewise-linear estimator, and $T_{0} = 0.5 T$ for the shrinkage estimator.} 
\end{figure}

The bias of this estimator is given by 
\begin{equation}
b_{n}^{\text{(pl)}} (x_{n}) =(1- \alpha) b_{n}^{\text{(ht)}}(x_{n})  ,
\end{equation}
and because $0\leq  (1-\alpha) \leq 1$, the bias and its derivative have smaller magnitude in comparison to hard thresholding, as can be verified from Figure~\ref{PL:BIASFIG}. The mean-square error is given by 
\begin{align}
\text{mse}_{n}^{\text{(pl)}}(x_{n}) = \sigma_{w}^{2} - & (1-\alpha) 2 x_{n} b_{n}^{\text{(ht)}}(x_{n}) \nonumber \\  - & (1-\alpha^{2}) \int_{-T}^{T} {y^{2} \, p_{w}(y-x_{n}) \, {\rm d}y} ,
\end{align}
and is plotted together with the mean-square error from hard thresholding in Figure~\ref{PL:MSEFIG}. It is clear from the plot that the piecewise-linear estimator has better worst-case mean-square error than hard thresholding, but this comes at the price of worse best-case mean-square error.

The second estimator we consider here is the ``semisoft shrinkage" estimator~\cite{GaoB:95} which replaces the discontinuities in the hard thresholding estimator with a linear segment connecting the left and right limit points, i.e., 
\begin{equation}
\hat{x}_{n}^{\text{(ss)}} (y_{n})
= \begin{cases} 0 , & \text{if } |y_{n}|< T_{0} ; \\
\beta (y_{n} - \text{sgn}(y_{n}) T_{0} ) & \text{if } T_{0} \leq |y| < T ; \\ 
y_{n} , & \text{if } |y_{n}|\geq T , \end{cases}
\label{SS:EST}
\end{equation}
where $\beta \equiv T/(T - T_{0}) > 1$ denotes the slope of the line segment shown in Figure~\ref{SS:ESTFIG}, and 
\begin{equation}
\text{sgn}(y) = \begin{cases} 1 , & y>0 \\ 0 , & y = 0 \\ -1 , & y<0 \end{cases} \label{signum}
\end{equation}
is the signum function. Note that the shrinkage estimator reduces to hard thresholding when $T = T_{0}$, but it does not otherwise encompass the piecewise-linear estimator in \eqref{PL:EST}.

The bias and mean-square error expressions for the shrinkage estimator are less tractable than the previous case. Nonetheless, its bias can be expressed as
\begin{align}
b_{n}^{\text{(ss)}} (x_{n}) &= b_{n}^{\text{(ht)}}(x_{n})  + \beta \int_{T_{0}}^{T} (y-T_{0}) \nonumber \\ & \times \bigl ( p_{w}(y-x_{n}) - p_{w}(y+x_{n}) \bigr ) \, {\rm d}y , 
\end{align}
and is plotted in Figure~\ref{PL:BIASFIG} together with those obtained from the estimators introduced thus far. The mean-square error of this shrinkage estimator takes on the form
\begin{equation}
 \text{mse}_{n}^{\text{(ss)}}(x_{n}) = \text{mse}_{n}^{\text{(ht)}}(x_{n}) + f(x_{n}) + f(-x_{n}) \label{SS:MSE}
\end{equation}
where 
\begin{equation}
f(x) \equiv \int_{T_{0}}^{T} \bigl( \beta^{2} (y-T_{0})^{2} - 2 x \beta (y - T_{0}) \bigr ) p_{w}(y-x) \, {\rm d}y. 
\end{equation}
This mean-square error is compared to that of the previous estimators in Figure~\ref{PL:MSEFIG}. It is seen that the shrinkage estimator has similar error to that obtained from \eqref{PL:EST}, but the peak of the oscillation is slightly skewed.

To demonstrate that piecewise-linear estimation improves performance over hard thresholding, we return to the vector-valued AWGN estimation problem and compare the mean-square error per symbol obtained with a hard thresholding estimator and with the two piecewise-linear estimators introduced above. We assume that there are $N$ basis coefficients to be estimated from the same number of observations. Furthermore, because both estimators have symmetric mean-square error as a function of the true value of the parameter, we shall assume with no loss of generality that $x_{n} \geq 0$ for all $n=1,\dots,N$. Finally, we assume that the true values of the coefficients---when sorted---have a decay rate governed by a generalized Gaussian function, i.e. we assume the coefficient sequence is given by 
\begin{equation}
x_{n} = \kappa(p) e^{-\left [\lambda (n-1) \right ]^{p}}
\end{equation}
for $n=1\dots,N$, where $p>0$ is the decay rate, $\kappa(p)>0$ is a scaling factor such that the energy of the sequence is equal for all values of $p$, and the closest integer to $1/\lambda$ is approximately the $e^{-1}$ attenuation point of the coefficients. The signal-to-noise ratio of such a sequence is defined as the ratio of the total signal energy to the total noise energy, i.e.
\begin{equation}
\text{SNR} = \frac{\sum_{n=1}^{N} |x_n|^2}{N \sigma_{w}^{2}} .
\end{equation}

Recall that the estimators act on each coefficient separately and the noise on each coefficient is independent. Therefore, for both estimators, the total mean-square error is equal to the sum of the mean-square error from each coefficient, i.e.
\begin{equation}
E[\lVert \xhat^{\text{(m)}} (\boldsymbol{y})-\boldsymbol{x} \rVert^{2}] = \sum_{n=1}^{N} E\Bigl [\bigl (\hat{x}_{n}^{\text{(m)}}(y_{n}) - x_{n} \bigr)^{2} \Bigr ] ,
\end{equation}
where method $\text{m}$ is $\text{ht}$ for hard thresholding, $\text{pl}$ for piecewise linear, or $\text{ss}$ for semisoft shrinkage. Furthermore, the average mean-square error per symbol is defined as
\begin{equation}
E \bigl [\lVert \xhat^{\text{(m)}} (\boldsymbol{y})-\boldsymbol{x} \rVert^{2} \bigr ]/N .
\end{equation}

\begin{figure}
\centering
\includegraphics[width=2.7in]{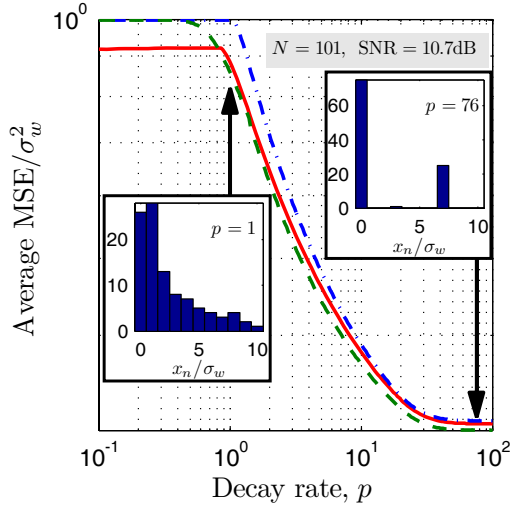}
\caption{Mean-square error per symbol of the three estimators analyzed herein. The dash-dotted (blue) curve is the mean-square error of hard thresholding estimators, the solid (red) curve is that of piecewise-linear estimators of form \eqref{PL:EST}, and the dashed (green) curve is the mean-square error of semisoft shrinkage estimators. Histograms refer to signal coefficients with decay rate $p=1$ and $p=75$ respectively. Optimization over threshold values and slope is performed numerically from analytic expressions (see footnote $4$ for details). Parameter values are $N=101$, $\lambda= 0.04$ and SNR=$10.7$dB.}
\label{PLHT:OPT}
\end{figure}

Figure~\ref{PLHT:OPT} plots the average mean-square error obtained from the three estimators as a function of the decay rate of coefficients when $N=101$, $\lambda=0.04$ and all the degrees of freedom for the estimators are optimized; i.e. the optimization is carried out over $T$ in \eqref{HT:EST}, $\alpha$ and $T$ in \eqref{PL:EST}, and $T_{0}$ and $T$ in \eqref{SS:EST}.\footnote{All coefficient sequences (each with different $p$) were normalized to the same energy (to attain identical $\text{SNR}$ in all sequences), where the normalization constant was chosen such that the largest coefficient over all $p$ was $10 \sigma_{w}$. Then, the optimization was performed numerically (using the analytical expressions given in Appendix~\ref{app2}), given the constraints $\alpha \in [0,1]$ and $T \geq T_{0} \geq 0$.} 

The plot shows that the average mean-square error improves for all values of the decay rate when either of the piecewise-linear estimators is utilized in place of the hard thresholding estimator. Let us consider the limiting cases. The histogram for $p=75$ shows that for fast decay rates the coefficients are clustered into two groups: a significant number of the coefficients are very close to $0$, while the remaining coefficients are grouped at a large value determined by the $\text{SNR}$ of the signal. Consequently, very few coefficients are in the intermediate region. Such a histogram is ideal for hard thresholding because the optimum threshold aligns the region with significant mean-square error in the gap between the two sets of coefficients. Thus, the error incurred per large coefficient becomes  $\sigma_{w}^{2}$  (the maximum-likelihood estimator limit), whereas coefficients that are approximately $0$ yield error that is a fraction of $\sigma_{w}^{2}$. In this regime, the optimal slope $\alpha$ in \eqref{PL:EST} approaches $0$, while the optimal threshold equals that of hard thresholding. Hence the mean-square error from the two estimators converge for $p \gg 1$. On the other hand, the smoother transition at the thresholding boundaries in the semisoft shrinkage estimator reduces the error incurred from the few coefficients that fall within the intermediate region, without significant impact on the errors incurred from the two main clusters of coefficients. Hence, for $p \gg 1$ the shrinkage estimator slightly outperforms the other two estimators.   

In the opposite limiting case, the slow decay rate implies that the coefficients will be more spread out over the parameter space, as evidenced by the histogram for $p=1$. Consequently, the error incurred from the coefficients with intermediate values becomes prohibitively large when either the hard thresholding or the shrinkage estimator is utilized. Hence, the optimal threshold parameters for both of these estimators, when $p \ll 1$, leads to the maximum-likelihood estimator (i.e. $T=0$ for hard thresholding and $T_{0}= 0$ for shrinkage). On the other hand, the optimal values for the piecewise-linear estimator in \eqref{PL:EST} turns out to be a large threshold value combined with a slope slightly smaller than unity. Therefore, an estimator that is linear over a significant portion of the coefficients, but with a more conservative slope than the maximum-likelihood estimator, reduces the average mean-square error in comparison to the that obtained from the maximum-likelihood estimator. This result is not entirely unexpected, as it has been shown previously that a (biased) linear estimator with slope less than unity performs better than the maximum-likelihood estimator over the entire parameter space~\cite{Eldar:06}. In the case considered herein, the slow decay rate implies that the coefficients will be spread out over the parameter space, thus the estimator that minimizes the average mean-square error per coefficient must perform well over a large subspace of the parameter space, and this is consistent with the optimization criterion considered in \cite{Eldar:06}. 

\begin{figure}
\centering
\subfigure[Slope]{
		\includegraphics[width=2.7in]{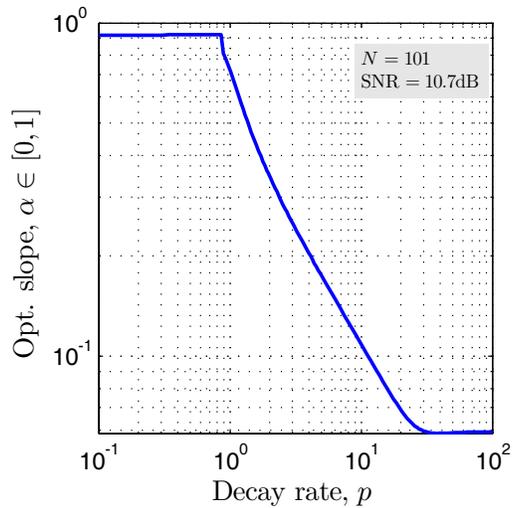}
		\label{PLHT:SLOPE}
	}\\
\subfigure[Thresholds]{
		\includegraphics[width=2.7in]{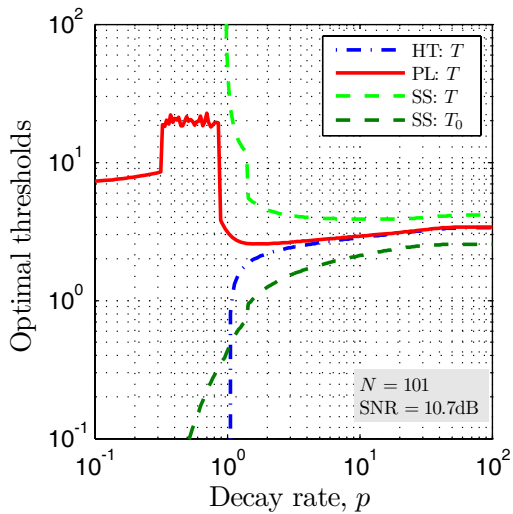}
		\label{PLHT:THRESH}
	}
	\caption{The optimum slope of the piecewise-linear estimator and the optimal threshold values of all estimators are plotted as a function of the decay rate of the signal coefficients. Values are determined numerically from analytical expressions given in Appendix~\ref{app2}. `HT' denotes the hard thresholding estimator, `PL' is the piecewise-linear estimator in \eqref{PL:EST}, and `SS' denotes the semisoft shrinkage estimator. Parameter values are $N=101$, $\lambda= 0.04$ and SNR=$10.7$dB.} 
\label{PLHT:OPTIMALS}
	\end{figure}

The optimal values of the degrees of freedom for all three estimators are plotted in Figure~\ref{PLHT:OPTIMALS} as a function of the decay rate $p$.

\section{Discussion}
\label{sec:discussion}

In this paper we have provided an estimation theoretic study of non-random signal estimation in order to deepen our understanding of some of the common results encountered in wavelet-based estimation techniques. We focused on the problem of  estimating the basis expansion coefficients of a signal when the coefficients are corrupted with additive white Gaussian noise. 

The main results developed in this paper can be summarized as follows. The mean-square error lower bound that applies to unbiased estimators indicates that linear estimation is optimal, and furthermore, the basis choice for decomposing the signal does not affect optimal performance. This is, of course, an expected result. Optimal processing of Gaussian random vectors is always linear. Furthermore, if one obtains an optimal solution in one basis, a non-singular transformation of the coordinate-system does not affect the minimum achievable mean-square error, because the optimal estimator in the new basis will simply invert the transformation and apply the former estimator to achieve the same minimum.

The results for biased estimators however differ notably from those for unbiased estimators. Our analysis of hard thresholding as an example of a biased estimator shows that biased estimators are not constrained by the unbiased version of the Cram\'{e}r-Rao bound. Furthermore, the extension of this bound for biased estimators does not yield achievable lower bounds on performance. Hence, optimality arguments for biased estimators are inevitably more heuristic. Our analysis for hard thresholding demonstrated that basis representation indeed does affect performance when such an estimator is used, as the basis coefficients must be well separated into values that are very close to zero (approximately within one standard-deviation of the noise) and values that are large (significantly larger than the noise standard deviation and the threshold value) to obtain mean-square error smaller than the error of a maximum-likelihood estimator. In other words, the decay rate of the sorted basis coefficients must be fast. Therefore, for the class of signals that have fast-decaying wavelet coefficients, wavelet-basis decompositions in conjunction with hard thresholding will be effective in denoising the observed signal.

Nevertheless, the Cram\'{e}r-Rao bound for biased estimators provides additional information on how the bias of an estimator affects mean-square error performance. In particular, through our analysis of this bound for the hard thresholding estimator we motivated piecewise-linear estimators and subsequently demonstrated that they achieve smaller average mean-square error when the decay rate of the coefficients are governed by a generalized Gaussian distribution. 

In summary, although prior literature provides abundant analysis on the reduction of mean-square error by utilizing wavelet basis expansions and thresholding estimators, additional insight can be obtained by connecting the recent advances in wavelet-based techniques with more traditional estimation theoretic analysis. In this paper we have provided such a connection through non-random signal estimation theory for the AWGN problem, and we have used our analysis to demonstrate that piecewise-linear estimators improve the average mean-square error attained via hard thresholding.

\appendices
\section{Analytical Expressions for Bias and Mean-Square Error in Hard Thresholding}
\label{app2}

Let us define 
\begin{align}
x_{S} &= (x_{n} + T)/(\sqrt{2}\, \sigma_{w}),\\
x_{D} &= (x_{n} - T)/(\sqrt{2}\, \sigma_{w}).
\end{align}
Evaluating the bias for the hard thresholding estimator from \eqref{HT:BIAS} gives
\begin{align}
b_{n}^{\text{(ht)}}(x_{n})/\sigma_{w} = (2\pi)^{-1/2} & \left(e^{-x_{D}^{2}} - e^{-x_{S}^{2}}\right) \nonumber \\ - & (x_{n}/\sigma_{w}) \bigl (Q(x_{D}) - Q(x_{S})  \bigr) ,
\end{align}
where
\begin{equation}
Q(x) \equiv \pi^{-1/2} \int_{x}^{\infty} e^{-t^2} {\rm d}t . 
\end{equation}

The mean-square error, on the other hand, is obtained from \eqref{HT:MSE} as
\begin{align}
\lefteqn{\text{mse}_{n}^{\text{(ht)}}(x_{n})/ \sigma_{w}^{2} = 1 + (x_{n}/\sigma_{w})^{2} \bigl (Q(x_{D}) - Q(x_{S}) \bigr)} \nonumber \\ & &+ \frac{1}{2} \Bigl (\text{sgn}(x_{D}) \Gamma_{\text{inc}}\bigl(x_{D}^{2},3/2 \bigr )\!-\! \text{sgn}(x_{S}) \Gamma_{\text{inc}}\bigl (x_{S}^{2},3/2\bigr) \Bigr), \label{HT:MSEEXP}
\end{align}
where $\text{sgn}(x)$ is defined in \eqref{signum}, and
\begin{equation}
\Gamma_{\text{inc}}\bigl (x, 3/2 \bigr) \equiv 2/\sqrt{\pi} \int_{0}^{x} t^{1/2} e^{-t} {\rm d}t
\end{equation}
is the Gamma distribution of order $3/2$.

To provide analytical expressions for the shrinkage estimator we must define two new dimensionless variables,
\begin{align}
\xi_{S} &= (x_{n} + T_{0})/(\sqrt{2}\, \sigma_{w}),\\
\xi_{D} &= (x_{n} - T_{0})/(\sqrt{2}\, \sigma_{w}).
\end{align}
Now the bias is given by 
\begin{align}
b_n^{\text{(ss)}} (x_{n}) & = \beta b_{n}^{\text{(ht)}} (x_{n} ; T_{0}) - (\beta - 1) b_{n}^{\text{(ht)}} (x_{n} ; T) \nonumber \\ & - \beta T_{0} \bigl ( Q(x_{D}) - Q(\xi_D) + Q(x_{S}) - Q(\xi_{S}) \bigr ),
\end{align}
where the second argument in $b_{n}^{\text{(ht)}}$ indicates the threshold value. The mean-square error expression \eqref{SS:MSE} depends on \eqref{HT:MSEEXP} and the function
\begin{align}
\lefteqn{f(x) / \beta^{2} \sigma_{w}^{2} = } \nonumber \\  &  \Bigl ( x_{D} (1- 2T_{0}/T) e^{-x_{D}^{2}} -  \bigl (\xi_{S} \!-\! \sqrt{2}T_{0} x/(\sigma_{w} T) \bigr ) e^{-\xi_{D}^{2}} \Bigr)/\sqrt{\pi}  \nonumber \\ &    + \Bigl (1 - 2 \xi_{D} \bigl (\xi_{S} \!-\! \sqrt{2} T_{0} x/(\sigma_{w} T) \bigr ) \Bigr ) \bigl ( Q(x_{D}) \!-\! Q(\xi_{D}) \bigr ).
\end{align}

\section{``Optimal Oracle" Bound}
\label{app3}
We simply reproduce the derivation in \cite[Ch.~10]{Mallat:99}. Consider a scalar estimator of the form
\begin{equation}
\hat{x}(y) = a y ,
\label{ORACLEEST}
\end{equation}
where $a \in \mathbb{R}$ is deterministic, for the scalar AWGN estimation problem. Then
\begin{equation}
E \bigl [ (\hat{x}(y) - x)^{2} \bigr ] = a^2 \sigma_{w}^{2} + (1-a)^2 x^{2} .
\end{equation}
Differentiating this expression with respect to $a$ and setting it to $0$, we find the value that minimizes the mean-square error as
\begin{equation}
a_{\text{opt}} = \frac{x^{2}}{x^{2} + \sigma_{w}^{2}} ,
\end{equation}
and the minimum mean-square error is
\begin{equation}
\min_{a\in\mathbb{R}} E [(a y -x)^2 ]  = \frac{\sigma_{w}^{2} x^{2}}{x^{2} + \sigma_{w}^{2}} .
\label{Oracle:MSE}
\end{equation}
Equation \eqref{Oracle:MSE} is the oracle bound used in \cite[Ch.~10]{Mallat:99} and plotted in Figure~\ref{HT:LBFIG}. Because $a_{\text{opt}}$ depends on $x$ (which is unknown), such an estimator is not feasible. Hence this mean-square error is not achievable by any feasible estimator of form given in \eqref{ORACLEEST}.

\bibliographystyle{IEEEtran} 
\bibliography{bibl}

\end{document}